\documentclass[10pt]{amsart}
\usepackage{amssymb,xypic} 

\newtheorem{thm}{Theorem}

\newtheorem{prop}[thm]{Proposition}
\newtheorem{lem}[thm]{Lemma}

\theoremstyle{definition}
\newtheorem{defn}[thm]{Definition}
\newtheorem{exam}[thm]{Examples}

\newcommand{\Gl}{\mathrm{GL}}

\newcommand{\g}{\mathfrak{g}}
\newcommand{\tr}{\mathrm{Tr}}
\newcommand{\Tr}{\mathrm{Tr}}

\newcommand{\MC}{\mathcal{MC}}
\newcommand{\mc}{\mathcal{MC}}
\renewcommand{\d}{\mathrm{d}}
\newcommand{\A}{A}
\newcommand{\B}{\mathcal B}
\newcommand{\R}{\mathbb R}
\newcommand{\C}{\mathbb C}
\newcommand{\Z}{\mathbb {Z}}

\newcommand{\db}{\overline\partial}

\title[Algebraic string bracket as a Poisson bracket]{Algebraic string bracket as a Poisson bracket}
\author{Hossein Abbaspour\and Thomas Tradler \and Mahmoud Zeinalian}

\keywords{Free loop space, Cyclic homology, Maurer Cartan, Symplectic reduction}

\begin{document} 

\bibliographystyle{h-elsevier2}
\begin{abstract}
In this paper we construct a Lie algebra representation of the algebraic string bracket on negative cyclic cohomology of an associative algebra with appropriate duality. This is a generalized algebraic version of the main theorem of [AZ] which extends Goldman's results using string topology operations.The main result can be applied to the de Rham complex of a smooth manifold as well as the Dolbeault resolution of the endomorphisms of a holomorphic bundle on a Calabi-Yau manifold.
\end{abstract}
\maketitle
\tableofcontents
\section{Introduction}

Goldman original work \cite{Gold} on the Lie algebra of free homotopy classes of oriented closed curves on an oriented surface was extensively generalized through the introduction of String Topology by Chas and Sullivan \cite{CS1}. In particular, they generalized this Lie bracket to one on the equivariant homology of the free loop space of a compact and oriented manifold $M$. From the beginning, it was clear that this bracket had a deep relation to the holonomy map on a vector bundle; see \cite{Gold, CFP, CCR, CR}. This relation was the subject of a paper, \cite{AZ}, by the first and third author. It was shown there that using Chen's iterated integral one obtains a map of Lie algebras from the equivariant homology of the free loop space to the space of functions on a space of generalized flat connections.

Algebraic analogues of string topology Lie algebra have also been considered in recent years. Jones \cite{J} had shown that for a simply connected topological space $X$ the equivariant homology of the free loop space is isomorphic to the negative cyclic cohomology of the algebra of cochains on $X$. Using this, and Connes long exact sequence relating negative cyclic cohomology and Hochschild cohomology, together with the BV algebra on Hochschild cohomology, Menichi \cite{Men} deduced a Lie bracket on the negative cyclic cohomology in a way similar to the one in string topology \cite[Section 6]{CS1}. 

The starting point for this work was to obtain a generalization of the the results in \cite{AZ} and place it in a more algebraic setting where the equivariant homology of the loop space is replaced by negative cyclic cohomology. A suitable setting for this is to consider a unital differential graded algebra $A$ over a field $k = \R ~\text{or}~ \C$, with a reasonable trace $\tr:A\to k$. Using the results of \cite{T}, the above assumptions imply an isomorphism of the Hochschild cohomologies of $A$ with values in $A$ and its dual $A^*$, $HH^\bullet(A,A)\cong HH^\bullet(A,A^*)$, such that the cup product on $HH^\bullet(A,A)$ and the dual of Connes $B$-operator on $HH^\bullet(A,A^*)$ make these spaces into a BV algebra. This BV algebra, together with a Connes long exact sequence between the Hochschild cohomology $HH^\bullet(A,A^*)$ and negative cyclic cohomology $HC_-^\bullet(A)$, imply a Lie algebra structure on $HC_-^\bullet(A)$ by a theorem of Menichi's \cite[Proposition 7.1]{Men}, which is based on a similar marking/erasing result of Chas and Sullivan \cite[Theorem 6.1]{CS1}.

Now, using work of Gan and Ginzburg in \cite{GG}, we may look at the moduli space of Maurer-Cartan solutions,
\begin{equation}\label{MC-intro}
 \MC= \{a\in A^{odd}~|~ da+a\cdot a=0\}/\sim
\end{equation}
Since we only consider odd elements, the trace induces a symplectic structure $\omega$ on $\MC$, and thus one can define a Poisson bracket on the function ring $\mathcal O(\MC)$ of $\MC$. More details of this construction will be given in Section \ref{MC-section}.

We may connect the two sides of the above discussion via a canonical map $\{a\in A^{odd}~|~da+a\cdot a=0\}\to HC^-_\bullet(A), a\mapsto \sum_{n\geq 0} 1\otimes a^{\otimes n}$, and dualizing this gives a map $\rho :HC_-^\bullet(\A) \to \mathcal O(\MC)$. We may now compare the two Lie algebras from above. Our main result then states, that the brackets are indeed preserved.
\begin{thm}\label{main-theorem}
$\rho: HC_-^{2\bullet}(\A) \to \mathcal O(\MC)$ is a map of Lie algebras.
\end{thm}
In a special case considered in \cite{AZ} this map becomes the generalized holonomy map from the equivariant homology of the free loop space of $M$ to the space of functions on the moduli space of generalized flat connections on a vector bundle $E\to M$. In fact one has a commutative diagram,  
\begin{equation}\label{diag1}
\xymatrix{HC_-^{2\bullet}(\A)  \ar[rr]^{\rho} & & \mathcal{O(\MC)}\\
& H_{2\bullet} ^{S^1}(LM) \ar[ur]_{\Psi} \ar[ul]^{\sigma\quad}& }
\end{equation}
where $\Psi$ is the generalized holonomy dicussed in \cite{AZ} and $\sigma$ comes from Chen's iterated integral map, as described in Section \ref{gen-hol}. In particular, for $\dim M=2$, this recovers Goldman's results on the space of flat connections on a surface. 

Another motivation of this work is to study string topology in a holomorphic setting via the moduli stack of the holomorphic structure on a fixed complex bundle $E\rightarrow M$, where $M$ is a complex manifold. Algebraically, this will correspond to the choice of the algebra $A=\Omega^{0,*}(M, End(E))$, with the Dolbeault differential $\db$.  This discussion, once done at the chain level, relates to the algebraic structure of the B-model.

Finally, we remark, that the above discussion generalizes in a straight forward way to the case of a cyclic $A_\infty$ algebra $A$. This will be the topic of the last Section \ref{A-infty-section}. In fact, by the same reasoning as above, we obtain the Lie bracket on the negative cyclic cohomology $HC^\bullet_-(A)$. Also, by symmetrization we may associate an $L_\infty$ algebra to $A$, which induces a Maurer-Cartan space similar to \eqref{MC-intro}.   We find, that the canonical map $\rho$ is still well-defined, such that Theorem \ref{main-theorem}  also remains valid in this generalized setting. 

\textbf{Notation:} For a map $F$ of complexes, $F_\bullet$ (resp. $F^\bullet$) denotes the induced map in homology (resp. cohomology). 

\textbf{Acknowledgments:} 
 The authors would like to  thank Victor Ginzburg and Luc Menichi for useful discussions and correspondence on this topic. The authors were partially supported by the Max-Planck Institute in Bonn and the second author warmly thanks the Laboratoire Jean Leray at the University of Nantes for their invitation throught the Matpyl program.

\section{The Lie algebra $HC^\bullet_-(\A)$}\label{sectionHC-}

In this section, we recall the Lie algebra structure of the negative cyclic cohomology $HC^\bullet_-(A)$, for a dga $(A,d,\cdot)$ with a trace $\tr:A\to k$. The Lie bracket comes from the long exact sequence that relates negative cyclic (co-)homology to Hochschild (co-)homology. For simplicity, we will work in the normalized setting.

\begin{defn}\label{Hoch}
Let $(A=\bigoplus_{i\in \Z}A^i,d:A^i\to A^{i+1},\cdot)$ be a differential graded associative algebra over a field $k$, and let $M=\bigoplus_{i\in \Z}M^i$ be a differential graded $A$-bimodule. The (normalized) Hochschild chain complex defined as,

\begin{equation}\label{C(A,M)}
 \bar C_\bullet(A,M):=\prod_{n\geq 0} M\otimes \bar A ^{\otimes n}, 
\end{equation}
where $\bar A=A/k$, and $s$ denotes shifting down by one. The boundary $\delta:\bar C_\bullet(A,M)\to \bar C_{\bullet+1}(A,M)$ is defined by,
\begin{multline*}
\delta(a_0\otimes a_1\otimes\cdots\otimes a_n):= \sum_{i=0}^n (-1)^{\epsilon_i} a_0\otimes\cdots\otimes d(a_i)\otimes \cdots\otimes a_n \\
+\sum_{i=0}^{n-1} (-1)^{\epsilon_i} a_0\otimes\cdots\otimes (a_i\cdot a_{i+1})\otimes \cdots\otimes a_n -(-1)^{\epsilon'_n} (a_n\cdot a_0)\otimes a_1\otimes\cdots\otimes a_{n-1},
\end{multline*}
where $a_0\in M$, $a_1,\cdots,a_n\in A, \epsilon_0=|a_0|, \epsilon_i=(|a_0|+\cdots+|a_{i-1}|+i-1)$, and $\epsilon'_n=\left(|a_n|+1\right)\cdot
\left(|a_0|+\cdots+|a_{n-1}|+n-1\right)$.
Note that the differential is well defined; see \cite{L}.
Similarly, the (normalized) Hochschild cochain complex is defined by,
\begin{equation}\label{C*(A,M)}
 \bar C^n(A,M):=\Big\{f:s\bar A^{\otimes n}\to M \quad\Big| \quad f(a_{1}\otimes \cdots\otimes a_i\otimes \cdots\otimes a_n)=0, \text{if } a_i=1 \Big\}, 
\end{equation}
where the differential $\delta^*:\bar C^\bullet(A,M)\to \bar C^{\bullet-1}(A,M)$ is given by,
\begin{multline*}
(\delta^* f) (a_1\otimes\cdots\otimes a_n):= \sum_{i=1}^n (-1)^{|f|+\epsilon_i} f(a_1\otimes\cdots \otimes da_i\otimes \cdots\otimes a_n) \\ 
+d ( f (a_1\otimes\cdots\otimes a_n))+ \sum_{i=1}^{n-1} (-1)^{|f|+\epsilon_i} f(a_1\otimes\cdots \otimes (a_i\cdot a_{i+1})\otimes \cdots\otimes a_n) \\
+(-1)^{|f|(|a_1|+1)} a_1\cdot f(a_1\otimes\cdots\otimes a_{n})+ (-1)^{|f|+\epsilon_n} f(a_1\otimes\cdots\otimes a_{n-1}) \cdot a_n.
\end{multline*}
The respective (co-)homology theories are denoted by
$$ HH_\bullet(A,M)=H(\bar C_\bullet(A,M),\delta),\quad HH^\bullet(A,M)=H(\bar C^\bullet(A,M),\delta^*). $$ Denoting by $A^*=Hom(A,k)$ the graded dual of $A$, we see that the dual of $\bar C_\bullet(A,A)$ is given by $ \bar C^\bullet(A,A^*)$. Recall furthermore, that there is a cup product $\cup$ on $\bar C^\bullet(A,A)$ defined by
$$ (f\cup g)(a_1\otimes \cdots\otimes a_{m+n}):=f(a_1\otimes\cdots\otimes a_m)\cdot g(a_{m+1}\otimes\cdots\otimes a_{m+n}). $$

Next, we define the (normalized) negative cyclic chains $\overline{CC}_\bullet^-(A)$ of $A$ to be  the vector space $\bar C_\bullet(A,A)[[u]]$, where $u$ is of degree $+2$, and with differential $\delta+uB$, where $B:\bar C_\bullet(A,A)\to \bar C_{\bullet-1}(A,A)$ is Connes operator,
\begin{multline}\label{B}
 B(a_0\otimes a_1\otimes\cdots\otimes a_n):=\sum_{i=0}^n (-1)^{\epsilon_i}1\otimes a_i\otimes\cdots\otimes a_n\otimes a_0\otimes \cdots\otimes a_{i-1},\\
\text{where } \epsilon_i=(|a_i|+\cdots+|a_n|+n-i+1)(|a_0|+\cdots+|a_{i-1}|+i-1).
\end{multline}
Thus, every element of $\overline{CC}_n^-(A)$ is an infinte sum $\sum_{i=0}^{\infty} a_i
u^i\in \bar C_\bullet(A,A)[[u]]$, where $a_i\in \bar C_{n-2i}(A,A)$, $\delta$ acts on $a_i\in  \bar
C_\bullet(A,A)$, and $uB$ acts as
\begin{equation}\label{CC-}
 \cdots \stackrel {uB} \longleftarrow \bar C_\bullet(A,A) \cdot u^2 \stackrel {uB} \longleftarrow \bar  C_\bullet(A,A)\cdot u\stackrel {uB} \longleftarrow \bar C_\bullet(A,A).
\end{equation}

Dually, define the (normalized) negative cyclic cochains $\overline{CC}^\bullet_-(A)$ of $A$ by taking $\overline{CC}^\bullet_-(A)= \bar C^\bullet(A,A^*)\otimes k[v,v^{-1}]/ vk[v]$, where $v$ is an element of degree $-2$. Explicitly, the degree $n$ part $\overline{CC}^n_-(A)$ is represented by finite sums $\sum_{i= 0}^k a_iv^{-i}$ where $a_i\in \bar C^{n-2i}(A,A^*)$. The differential is given by $\delta^*+vB^*$, where $\delta^*$ acts on $\bar C^\bullet(A,A^*)$, and $vB^*$ acts as follows.
$$ \cdots \stackrel {v B^*} \longrightarrow \bar C^\bullet(A,A^*)\cdot v^{-2} \stackrel {v B^*} \longrightarrow \bar C^\bullet(A,A^*)\cdot v^{-1} \stackrel {v B^*} \longrightarrow \bar C^\bullet(A,A^*). $$
Note, that if $C_\bullet(A,A)$ is finite dimensional in each degree, then the graded dual of $\overline{CC}^n_-(A)$ is isomorphic to the chain complex $\overline{CC}_n^-(A)=Hom(\overline{CC}^n_-(A),k)$, see also \cite[Lemma 3.7]{HL}.
It is easy to see that $B^2= \delta B+B\delta=0$, and we define the associated (co-)homology theories by,
$$
HC^-_\bullet(A)=H(\overline{CC}_\bullet^-(A), \delta+u B), \quad
HC_-^\bullet(A)=H(\overline{CC}^\bullet_-(A), \delta^*+v B^*).
$$
\end{defn}
\begin{lem}\label{[u]vs[[u]]}
If $H_\bullet(A,A)$ is bounded from below, then both $\bar C_\bullet(A,A)[u]$ and $\bar C_\bullet(A,A)[[u]]$ with differential $\delta+uB$ calculate negative cyclic homology $HC^-_\bullet(A)$. 
\end{lem}
This lemma follows from a spectral sequence argument for the inclusion $\bar C_\bullet(A,A)[u]\hookrightarrow\bar C_\bullet(A,A)[[u]]$, similarly to \cite[Lemma 3.6]{HL}. Note, that our sign convention is opposite to the one from \cite{HL}, but in agreement with \cite{GJP}, since our differential $\delta:\bar C_\bullet(A,A)\to \bar C_{\bullet+1}(A,A)$ is of degree $+1$.

From now on, we additionally assume, that we also have a suitable trace map.
\begin{defn}\label{trace}
Let $\tr:A\to k$ be a trace map, satisfying $\tr(da)=0$ and $\tr(ab)=-(-1)^{|a|\cdot |b|}\tr (ba)$, for all $a,b\in A$. Assume furthermore that the map $\omega:A\to A^*$, $\omega(a)(b):=\tr (ab)$ is a bimodule map, which induces an isomorphism on homology $H(A)\to H(A^*)$. By abuse of language, we will also view $\omega$ as a map $\omega:A\otimes A\to k, \omega(a, b)= \tr(ab)$. In this case, $A$ is also called a \emph{symmetric algebra}.

Notice that $\omega:A\to A^*$ induces a morphism of the Hochschild complexes $\omega_\sharp:\bar C^\bullet(A,A)\to \bar C^\bullet(A,A^*)$ via composition $\omega_\sharp (f):= \omega\circ f$, which is an isomorphism on homology $\omega^\bullet_\sharp:H^\bullet(A,A)\to H^\bullet(A,A^*)$. We may thus transfer the cup product $\cup$ on $H^\bullet(A,A)$ to a product $\sqcup$ on $HH^\bullet(A,A^*)$, by setting $f\sqcup g:=\omega^\bullet_\sharp((\omega^\bullet_\sharp)^{-1} f\cup (\omega^\bullet_\sharp)^{-1} g)$. Define furthermore the operator $\Delta:HH^\bullet(A,A^*)\to HH^\bullet(A,A^*)$ as the dual of $B$ on homology.
Then we assume, that $(HH^\bullet(A,A^*),\sqcup, \Delta)$ is a BV-algebra, {\it i.e.} $\sqcup$ is a graded associative, commutative product, $\Delta^2=0$, and the bracket $\{a,b\}:=(-1)^{|a|}\Delta(a\sqcup b)-(-1)^{|a|}\Delta(a)\sqcup b-a\sqcup \Delta(b)$ is a derivation in each variable.

Recall from Menichi \cite{Men} that this BV-algebra induces a Lie algebra on the negative cyclic cohomology $HC^\bullet_-(\A)$ using the long exact sequences of Hochschild and negative cyclic cohomology. The inclusion $\overline{CC}^-_\bullet(\A) \stackrel {\times u}\to \overline{CC}^-_\bullet(\A)$ given by multiplication by $u$ has cokernel $\bar C_\bullet(\A,\A)$. We thus obtain a short exact sequence
 \begin{equation}\label{short1} 0\to \overline{CC}^-_\bullet(\A)\overset{\times u}{\longrightarrow} \overline{CC}^-_\bullet(\A)\to \bar C_\bullet(\A,\A)\to 0,
\end{equation}
which induces Connes long exact sequence of homology groups.
\begin{equation}\label{CC_HH}
 \cdots\to HH_n(\A,\A)\stackrel {\B_\bullet} \to HC^-_{n-1} (\A)\stackrel {} \to HC^-_{n+1} (\A)\stackrel{I_\bullet} \to HH_{n+1}(\A,\A) \stackrel {\B_\bullet} \to \cdots.
\end{equation}
Here, the projection to the $u^0$ term $I:\overline{CC}_\bullet^-(\A)\to \bar C_\bullet(A,A)$ induces the map $I_\bullet$, and the connecting map $\B_\bullet$,  is induced by the composition $\bar C_\bullet(A,A)\stackrel B \to \bar C_\bullet(A,A)\stackrel {inc} \to \overline{CC}^-_\bullet(A)$. Note, that unlike $inc\circ B:\bar C_\bullet(A,A)\to  \overline{CC}^-_\bullet(A)$, the inclusion $inc:\bar C_\bullet(A,A)\to \overline{CC}^-_\bullet(A)$ is not a chain map.

Dually, we have the short exact sequence
$$ 0\to \bar C^\bullet(\A,\A) \to \overline {CC}_-^\bullet(\A)\to \overline {CC}_-^\bullet(\A)\to 0, $$
inducing Connes long exact sequence of cohomology groups
\begin{equation}\label{CC^HH}
\cdots\to HH^n(\A,\A^*)\stackrel {I^\bullet} \to HC_-^{n} (\A)\stackrel {} \to HC_-^{n-2} (\A)\stackrel {\B^\bullet} \to HH^{n-1}(\A,\A^*) \stackrel {I^\bullet} \to \cdots.  
\end{equation}
Notice that the composition 
\begin{equation}\label{equBB}
B^\bullet=\B^\bullet\circ I^\bullet
\end{equation}
 is exactly the $\Delta$ operator of our BV-algebra on $HH^\bullet(A,A^*)$, so that we may obtain an induced Lie algebra from \cite[Lemma 7.2]{Men}, much like the marking/erasing situation in \cite{CS1}.
\begin{prop}[L. Menichi \cite{Men}]\label{Lie}
The bracket $\{a,b\}:=I^\bullet(\B^\bullet(a)\sqcup \B^\bullet(b))$ induces a Lie algebra structure on $HC^\bullet_-(\A)$.
\end{prop}
\end{defn}

We end this section with some examples of the above definitions.
\begin{exam} \label{exam1}
Let $M$ be a smooth, compact and oriented Riemannian manifold.
\begin{itemize}
\item
A first example is obtained by taking $A=\Omega^\bullet(M)$ the De Rham forms on $M$, $d=d_{DR}$ the exterior derivative on $A$, and $\tr (a):=\int_M a$.
\item
More generally, if $E \to M$ is a finite dimensional complex vector bundle over $M$, with a flat connection $\nabla$, then we may take $A=\Omega^\bullet(M,  End(E))$ with the usual differential $d_\nabla$.  Similarly, the trace is given by a combination of integration and trace in $End(E)$. The cyclic property of the trace guarantees that this induces an injective bimodule map $\omega:A\to A^*$ that is a quasi-isomorphism.
\item
Both of the above examples are special cases of ellitpic Calabi-Yau space as defined in \cite{Cos}. By definition, this means that we have a bundle of finite dimensional associative $\C$ algebras over $M$, whose algebra of sections is denoted by $A$. Furthermore, there is a differential operator $d:A\to A$, which is an odd derivative with $d^2=0$ making $A$ into an ellitpic complex, a $\C$ linear trace $\tr:A\to \C$, a hermitian metric $ A\otimes A \to \mathbb C$, and a complex antilinear, $C^\infty (M, \mathbb R)$ linear operator $\ast: A \to A$, satisfying certain natural conditions. It can be seen that this example satisfies the above assumptions. The details and other examples of elliptic Calabi-Yau spaces can be found in \cite{Cos} and \cite{DonTho98}.
\end{itemize}
\end{exam}

\section{Maurer-Cartan solutions}\label{MC-section}
In this section we  define the moduli space of Maurer Cartan solutions for a symmetric algebra $\A=\oplus_{i\geq 0} A^i$, and then explain its symplectic nature. The main reference for this section is the paper \cite{GG} by Gan-Ginzburg, together with Section 4 of \cite{AZ}. Let us assume $k=\R$ or $\C$.

For $a,b \in \A$ define the Lie bracket $[a,b]:=a\cdot b-(-1)^{|a|\cdot |b|}b\cdot a$ and the bilinear form $\omega(a,b):= \Tr(ab)$. The first remark is that $(\A=\A^{odd}\oplus \A^{even}, d, [\cdot,\cdot], \omega)$ is a \emph{cyclic differential graded Lie algebra} as it is defined in Section 4 of \cite{AZ}, therefore all results in \cite{GG} applies here to define the Maurer-Cartan solutions. 

\begin{defn}
We define the Maurer-Cartan moduli stack as 
\begin{eqnarray*}
 MC&:=& \{a\in \A^{odd}~|~ da+\frac{1}{2}[a,a]=da +a\cdot a =0\}, \text{ and} \\
 \MC&:=& MC\big/\sim, 
\end{eqnarray*} 
where the equivalence is generated by the infinitesimal action of $\A^{0}$ on $\A$, where for $a\in A^{0}$, the vector field $\xi_{x}$ on $\A$ is defined by,
$$
\xi_{x}(a)=[x,a]-dx.
$$
\end{defn}
Recall that $\omega$ is a symplectic form and the infinitesimal action is Hamiltonian. Moreover, the map $\mu: a\mapsto \phi_{a}\in (\A^{even})^*$,  where
$$
 \phi(x)=\omega(da+\frac{1}{2}[a,a],x),
$$
is the moment map corresponding to the Hamiltonian action above. 
One should think of the tangent space $T_{[a]}\MC$ at a class $[a]$ as the 3-term complex
\begin{equation}\label{3term}
T_{[a]}\MC: \quad  T_{[a]}^{-1} \MC:=\A^{even}\overset{\xi(a)}{\longrightarrow}T_{[a]}^{0} \MC:=T_{a} A^{odd}=A^{odd}\overset{\mu'_{a}}{\longrightarrow} T_{[a]}^{1} \MC:={\A^{even}}^*,
\end{equation}
graded by -1, 0 and 1. Here $\xi(a)$ is the map $x\mapsto \xi_{x}(a)$. The $\ker \mu'$ is the Zarisky tangent space to $MC$ and the image of $\xi(a)$ accounts for the tangent space of the action orbit. Ideally, when $0$ is a regular value for $\mu$ and the infinitesimal action of $\A^{even}$ on $MC=\mu^{-1}(0)$ is free, this compex is concentrated in degree zero and the Zarisky tangent space to $\mc$ at $[a]$ is the cohomology group $H^0(T_{[a]}\mc)=H^*(\A^{odd}, d_{a})$ where $d_{a}b=db+[a,b]$. 

Note that  3-term complex (\ref{3term}) is self-dual where the self-duality at the middle term is given by the symplectic form
\begin{equation}\label{2form}
\omega(X_a,Y_a):=  \tr (X_a \cdot Y_a) \in k.
\end{equation}
By assumption from the previous section, $\omega$ is non-degenerate.
This gives rise to an isomorphism $ T_{[a]}\MC \overset{\simeq}{\rightarrow} (T_{[a]}\MC)^*$ and equips  $ (T_{[a]}\MC)$ with a symplectic form given by (\ref{2form}). In the case of a nonsingular point $[a]$ this is the usual pairing on $H^0(T_{[a]}\mc)=H(\A^{odd}, d_{a})$ induced by $\omega$.

The function space $\mathcal O(\MC)$ is defined to be the subspace of $\mathcal{O} (MC)$ invariant by the infinitesimal action. The symplectic form allows us to define the Hamiltonian vector field $X^\psi$ of a function $\psi\in \mathcal O(\MC)$ via
$$ \omega(X_a^\psi,Y_a)=\d\psi_a(Y_a):=\underset{t\rightarrow 0}{\lim}\frac{d}{dt}\psi(a+tY_a) ,\quad\quad \forall\,\, Y_a\in T_{[a]}^{1}\MC. $$
We then define the Poisson bracket on $\mathcal O(\MC)$ by,
$$ \{\psi,\chi\}:= \omega(X^\psi,X^\chi)= \tr(X^\psi \cdot X^\chi). $$

\section{The induced Lie map} \label{sec-map}

In this section, we define a map $\rho:HC_-^{2\bullet}(\A)\to \mathcal O(\MC)$, and prove it respects the brackets. We start by defining a map $P:MC\to \bar C_\bullet(\A,\A)$, and in turn the map $R:MC\to \overline{CC}_\bullet ^-(\A)$ which factors through $P$. Dualizing $R$ will induce the wanted map $\rho$.

\begin{defn}\label{PRrho}
Recall that $MC=\{a\in\A^{odd}~|~da+a\cdot a=0\}$ and $ \bar C_\bullet(\A,\A)=\prod_{n\geq 0} \A\otimes \bar\A^{\otimes n}$. Then, let $P:MC\to \bar C_\bullet(\A,\A)$ be given by the expression,
$$
P(a):=\sum_{i\geq 0}1\otimes a^{\otimes i}=(1\otimes 1)+(1\otimes a)+(1\otimes a\otimes a)+\cdots.
$$
Notice that for $a\in MC$, it is $\delta (P(a))=\sum 1\otimes a\otimes \cdots\otimes da\otimes \cdots \otimes a+ \sum 1\otimes a\otimes \cdots\otimes (a\cdot a)\otimes \cdots \otimes a=0$, due to the relation  $da+a\cdot a=0$ in $MC$. Thus, we obtain in fact a Hochschild homology class $[P(a)]\in HH_\bullet(\A,\A)$. 

Next, define the map $R:=inc\circ P$ as the composition $R:MC\stackrel P \to \bar C_\bullet(A,A)\stackrel {inc}\to \overline{CC}^-_\bullet(\A)$. Just as above, we have that $\delta(R(a))=0$, and since we are in the normalized setting, we see that $B(R(a))=0$, so that $(\delta+uB)(R(a))=0$. The induced negative cyclic homology class is again denoted by $[R(a)]\in HC^-_\bullet(\A)$. It is immediate to see that under the long exact sequence \eqref{CC_HH}, we have that $I(R(a))=P(a)$.

Using the pairing between between negative cyclic homology and negative cyclic cohomology, $\langle \cdot,\cdot \rangle:HC_-^\bullet (\A)\otimes HC^-_\bullet(\A)\to k$, we define the map $\rho$ by
\begin{eqnarray*}
\rho:HC_-^\bullet (\A)&\to& \mathcal O(\MC),\\
\rho([\alpha])([a])&:= &\langle [\alpha], [R(a)] \rangle=\langle \alpha, R(a) \rangle, \quad \quad \text {for } [\alpha]\in HC_-^\bullet (\A), [a]\in\MC.
\end{eqnarray*}
To simplify notation, we will also write  $\rho(\alpha)$ instead of $\rho([\alpha])$.
\end{defn}
\begin{lem}\label{zzz}
$\rho$ is well-defined.
\end{lem}
\begin{proof}
We need to show that the value $\rho([\alpha])([a])=\langle\alpha,R(a)\rangle$ is independent of the choice of the representative $[a]\in\{x\in\A^{odd}~|~dx+x\cdot x=0\}/\sim$. Infinitesimally, this amounts to showing that $L_{X(b)}\rho(\alpha)(a)=0$, where $L_{X(b)}$ is the Lie derivative along a vector field in the direction $X(b)_a=db+[a, b]\in T_{[a]}MC$, for any $b\in \A^{even}$. To see this, note that
\begin{eqnarray*}
L_{X(b)}\rho(\alpha)(a) &=& (i_{X(b)}\circ d + d\circ i_{X(b)})\rho(\alpha)(a)\\
&=& i_{X(b)}\circ d (\rho(\alpha))(a)\\
&=&\langle\alpha,\frac {d}{dt}\Big|_{t=0} R(a+tX(b)_a)\rangle
\end{eqnarray*}
Now, for any $Y_a\in T_{[a]}MC$, we have 

\begin{multline}\label{d/dt R}
 \quad\quad\quad \frac{d}{dt}\Big|_{t=0} R(a+tY_a)=1\otimes Y_a+1\otimes Y_a\otimes a+1\otimes a\otimes Y_a+\cdots\\
 =\B (Y_a+(Y_a\otimes a)+ (Y_a\otimes a\otimes a)+\cdots),\quad\quad\quad\quad
 \end{multline}
where we used Connes operator $\B:\bar C_\bullet(A,A) \to \overline{CC}^-_\bullet(A)$ from in the long exact sequence \eqref{CC_HH} applied to $Y_a+(Y_a\otimes a)+(Y_a\otimes a\otimes a)\in \bar C_\bullet(\A,\A)$. Thus, setting $Y_a=X(b)_a=db+[a, b]$ in the above expression, we obtain
\begin{eqnarray*}
L_{X(b)}\rho(\alpha)(a) &=& \langle\alpha, \B \big( db + [a, b] + db\otimes a+[a, b] \otimes a\\&& \quad \quad\quad\quad +db\otimes a\otimes a+[a, b]\otimes a\otimes a+\cdots\big)\rangle\\
&=&\langle \alpha, \B\circ \delta(b+(b\otimes a)+(b\otimes a\otimes a)+\cdots)\rangle\\
&=&\langle \alpha,  \delta \circ \B(b+(b\otimes a)+(b\otimes a\otimes a)+\cdots)\rangle\\
&=&\langle \delta^* \alpha,  \B(b+(b\otimes a)+(b\otimes a\otimes a)+\cdots)\rangle\\  
&=&0.
\end{eqnarray*}
\end{proof}
We are now ready to prove our main theorem.
\setcounter{thm}{0}
\begin{thm}
$\rho:HC_-^{2\bullet}(\A) \to \mathcal O(\MC)$ is a map of Lie algebras.
\end{thm} 
\setcounter{thm}{9}
\begin{proof}
We saw in \eqref{d/dt R} that $\frac{d}{dt}\big|_{t=0} R(a+tY_a)=\B (Y_a+(Y_a\otimes a)+(Y_a\otimes a\otimes a)+\cdots)\in \overline{CC}^-_\bullet (\A)$, where $(Y_a+(Y_a\otimes a)+(Y_a\otimes a\otimes a)+\cdots)\in \bar C_\bullet(\A,\A)$ for $Y_a\in T_{[a]}\MC$. Therefore,
\begin{eqnarray*}
(d\rho(\alpha))_a(Y_a)&=&\langle\alpha,\frac{d}{dt}\Big|_{t=0} R(a+tY_a)\rangle\\ &=&\langle\alpha, \B (Y_a+(Y_a\otimes a)+(Y_a\otimes a\otimes a)+\cdots)\rangle\\
&=&\langle \B^*\alpha,Y_a+(Y_a\otimes a)+(Y_a\otimes a\otimes a)+\cdots\rangle\\ &=&(\B^*\alpha)(1+a+a\otimes a+\cdots)(Y_a),
\end{eqnarray*}
where $\alpha\in \overline{CC}_-^\bullet(A), \B^*\alpha\in \bar C^\bullet(\A,\A^*)$, and thus $(\B^*\alpha)(\sum_{i\geq 0} a^{\otimes i})\in \A^*$. Now, using the isomorphism $\omega_\sharp^\bullet:HH^\bullet(A,A)\to HH^\bullet(A,A^*)$ from definition \ref{trace}, we apply its inverse to obtain an element $[f_\alpha]:=(\omega_\sharp^\bullet)^{-1} \B^\bullet [\alpha] \in HH^\bullet(A,A)$. We then claim that the Hamiltonian vector field $X^{\rho(\alpha)}_a$ may be expressed as
\begin{equation}\label{aa}
 X^{\rho(\alpha)}_a=f_\alpha \Big(\sum_{i\geq 0} a^{\otimes i}\Big)\quad \in T_{[a]}\MC.
\end{equation}
This should be compared with \cite[Lemma 7.2]{AZ} and \cite[Proposition 3.7]{Gold}.
To this end, first note, that the relation $0=(\delta^* f)(\sum_{i\geq 0} a^{\otimes i})=d_a(f(\sum_{i\geq 0} a^{\otimes i}))$, for $f\in \bar C^\bullet(A,A)$, shows that $X^{\rho(\alpha)}_a$ given by equation \eqref{aa}, represents a well-defined class in $T_{[a]}\MC$. We show \eqref{aa}, by applying the non-degeneracy of $\omega$ in the following equation, which is valid for any $Y_a\in T_{[a]}\MC$,
\begin{multline*}
\omega(f_\alpha (\sum a^{\otimes i}), Y_a)
=\tr(f_\alpha (\sum a^{\otimes i})\cdot Y_a)
=(\omega_\sharp f_\alpha) (\sum a^{\otimes i})(Y_a) \\
=(\B^*\alpha)(\sum a^{\otimes i})(Y_a)=(d\rho(\alpha))_a(Y_a)=\omega(X_a^{\rho(\alpha)},Y_a).
\end{multline*}

Now, calculating the Lie bracket gives
\begin{eqnarray*}
\rho(\{\alpha,\beta\})(a)&=& \langle\{[\alpha],[\beta]\},[R(a)]\rangle\\
&=& \langle I^\bullet (\B^\bullet[\alpha] \sqcup \B^\bullet[\beta]),[R(a)]\rangle\\ 
&=& \langle I^\bullet \omega_\sharp^\bullet ((\omega_\sharp^\bullet)^{-1}\B^\bullet[\alpha] \cup (\omega_\sharp^\bullet)^{-1} \B^\bullet[\beta]),[R(a)]\rangle\\ 
&=& \langle \omega_\sharp^\bullet ([f_\alpha] \cup [f_\beta]), I_\bullet [R(a)]\rangle \\
&=& \langle \omega_\sharp^\bullet ([f_\alpha] \cup [f_\beta]), [P(a)]\rangle.
\end{eqnarray*}
To evaluate this expression, note that for $f_\alpha:\bar A^{\otimes m}\to \A$ and $f_\beta:\bar A^{\otimes n}\to \A$, $\omega_\sharp^\bullet ([f_\alpha] \cup [f_\beta])$ is represented by the composition
$$ \bar A^{\otimes m+n}\stackrel {f_\alpha\otimes f_\beta}\longrightarrow \A\otimes \A\stackrel {\cdot} \to A\stackrel {\omega} \to A^*. $$
The first arrow with $f_\alpha\otimes f_\beta$ applied to $P(a)=1+(1\otimes a)+(1\otimes a\otimes a)+\cdots\in\prod_{i\geq 0} \A\otimes \bar\A^{\otimes i}$ then gives an expression, where we apply $a$ to all possible inputs in $\bar A^{\otimes n+m}$. To this, we then apply the product in $A$, and apply $\omega$ with input $1\in\A$, since $P(a)=1\otimes(\cdots)$. We thus obtain
\begin{eqnarray*}
\rho(\{\alpha,\beta\})(a)&=& \tr \Big(f_\alpha (1+a+a\otimes a+\cdots)\cdot f_\beta (1+a+a\otimes a+\cdots)\cdot 1 \Big) \\
&\stackrel {\eqref{aa}} {=}& \tr(X_a^{\rho(\alpha)} \cdot X_a^{\rho(\beta)}) = \omega(X_a^{\rho(\alpha)}, X_a^{\rho(\beta)}) =\{\rho(\alpha),\rho(\beta)\} (a).
\end{eqnarray*}
This is the claim of the theorem.
\end{proof}

\section{Comparison with generalized holonomy} \label{gen-hol}
In this section we compare the map $\rho$ with the generalized holonomy map $\Psi$ studied in \cite{AZ}. The relationship may be summarized in the diagram \eqref{diag1}. This shows how a special case the result of this paper relates to the main theorem of \cite{AZ}. The map $\Tr:A\rightarrow \C$ is  induced by the trace function on $\g \subseteq \Gl (n,\C)$ and integration of forms on $M$; see Example \ref{exam1}.

Our model of $S^1$-equivariant de Rham forms of $LM$ is $(\Omega(LM)[u], d+u\Delta )$ where  $u$ is a generator of degree 2 and  $\Delta: \Omega^\bullet(LM) \rightarrow \Omega^{\bullet-1}(LM)$ is the map induced by the $S^1$ action on $LM$; see \cite{GJP}. This model is quasi-isomorphic to the small Cartan model  $(\Omega_{inv}(LM)[u], d+i_X u)$ for the $S^1$ action, where $X$ is the fundamental vector field generated by the natural action of $S^1$.  The quasi-isomorphism is given by the averaging map $\Omega^\bullet(LM)\rightarrow \Omega_{inv}^\bullet(LM)$. More explicitly, for $\omega \in \Omega^\bullet(LM)$, $\Delta(\omega)$ is given by,
$$
\Delta(\omega)=\int _{fibre} ev^*(\omega)\in \Omega^{\bullet-1}(LM)
$$
\begin{equation}
\xymatrix{S^1\times LM \ar[r]^{ev} \ar[d]^{\pi} & LM \\ LM & } 
\end{equation}

Chen's iterated integral map and the trace map on $\g$ (see (6.3) \cite{AZ}, and Theorem A in \cite{GJP}) yields a map, which we denote by,
 $$
 S: (\bar C_\bullet(A,A), \delta) \rightarrow (\Omega^\bullet(LM), d).
 $$
$S$ induces the map $S^{HH}:HH_\bullet(\A,\A)\longrightarrow  H^\bullet (LM)$ on homology, and, after applying the pairing between homology and cohomology groups, we get,
$$
H_\bullet (LM)  \overset{\sigma^{HH}}{\longrightarrow}  HH^\bullet(\A,\A^*).
$$
Extending $S$ by $u$-linearity, we obtain a map, which we denote by abuse of notation by the same letter,
 $$
 S: (\bar C_\bullet(A,A)[u], \delta+uB) \rightarrow (\Omega^\bullet(LM)[u],d+u\Delta ).
 $$
Since, by Lemma \ref{[u]vs[[u]]}, $ (\bar C_\bullet(A,A)[u], \delta+uB )$ and $(\bar C_\bullet(A,A)[[u]], \delta+uB )$ are quasi-isomorphic in our setting, we obtain the induced map $S^{HC}:HC^-_\bullet(\A)\longrightarrow  H^\bullet_{S^1}(LM)$ on homology.
Composing $S^{HC}$ with the map $R:MC\to \overline{CC}_\bullet ^-(\A)=\bar C_\bullet(A,A)[u]$ from Section \ref{sec-map}, we get,
$$
MC \overset{R}{\longrightarrow}  HC^-_\bullet(\A)\overset{S^{HC}}{\longrightarrow}  H^\bullet_{S^1}(LM).
$$
Thus by duality, and using Lemma \ref{zzz}, we have,
$$
 H_\bullet ^{S^1}(LM)  \overset{\sigma=\sigma^{HC}}{\longrightarrow}  HC_-^\bullet(\A)\overset{\rho}{\longrightarrow}  \mathcal{O}(\mc).$$
The composition $\rho \circ \sigma $ is the generalized holonomy map $\Psi$ discussed in \cite {AZ}.
\begin{equation}
\xymatrix{HC_-^\bullet(\A)  \ar[rr]^{\rho} & & \mathcal{O(\MC)}\\
& H_\bullet^{S^1}(LM) \ar[ur]_{\Psi} \ar[ul]^{\sigma}& }
\end{equation}
It was proved in \cite{AZ}, that $\Psi$ is the morphism of Lie algebras. We will shortly see how this is  a consequence of Theorem \ref{main-theorem}. We first recall the following theorem.
\begin{thm}[S. Merkulov \cite{Mer}]\label{merkulov} The Chen integral induces a map of algebras $(H_\bullet(LM),\bullet)\rightarrow (HH^\bullet(A,A),\cup)$.
\end{thm}
Thus, by definition, $\sigma^{HH}:(H_\bullet(LM),\bullet)\rightarrow (HH^\bullet(A,A^*),\sqcup)$ is also a map of algebras. With this, we can now prove the following statement.
\begin{thm}  The map induced by the Chen iterated integrals $ \sigma: (H_\bullet^{S^1}(LM), \{\cdot,\cdot\})  \to  (HC_-^\bullet(\A), \{\cdot,\cdot\})$ is a map of Lie algebras. Here, the first bracket is the string bracket and the second one is defined in the statement of Proposition \ref{Lie}.
\end{thm}
\begin{proof}
The brackets on $H_\bullet^{S^1}(LM)$ and $HC_-^\bullet(A)$ are determined by the products on $H_\bullet(LM)$ and $HC^\bullet(A,A^*)$, together with the maps in the corrsponding Gysin long exact sequences. By Theorem \ref{merkulov}, it thus remains to show that the long exact sequences correspond to each other, {\it i.e.} that the following diagrams commute,
\begin{equation*}
\xymatrix{\cdots \ar[r] & H_\bullet^{S^1}(LM)\ar[d]^\sigma\ar[r]^ {m_\bullet} & H_{\bullet+1}(LM)\ar[d]^{\sigma^{HH}} \ar[r]^{e_\bullet}\ar[d] & H_{\bullet+1}^{S^1}(LM)\ar[d]^\sigma\ar[r]&H_{\bullet-1}^{S^1}(LM) \ar[d]^\sigma\ar[r]& \cdots\\
\cdots \ar[r] & HC^\bullet_-(A)\ar[r]^ { \mathcal B^\bullet} & HH^{\bullet+1}(A,A^*) \ar[r]^{I^\bullet} & HC^{\bullet+1}_-(A)\ar[r]&HC^{\bullet-1}_-(A) \ar[r]& \cdots}
\end{equation*}
Equivalently, we need to show the commutativity of the following dual sequence,
\begin{equation*}
\xymatrix{\cdots  \ar[r] & HC_\bullet^-(A)\ar[r]^{I_\bullet}\ar[d]^{S^{HC}} & HH_{\bullet}(A,A) \ar[r]^{\mathcal B_\bullet} \ar[d]^{S^{HC}}& HC_{\bullet-1}^-(A)\ar[r]\ar[d]^{S^{HH}} &HC_{\bullet-1}^-(A) \ar[r]\ar[d]^{S^{HC}}& \cdots  \\
\cdots  \ar[r] & H^\bullet_{S^1}(LM)\ar[r] ^{e^\bullet}& H^{\bullet}(LM) \ar[r]^{m^\bullet}& H^{\bullet-1}_{S^1}(LM)\ar[r]&H^{\bullet-1}_{S^1}(LM)\ar[r]& \cdots}
\end{equation*}
The top long exact sequence is induced by the short exact sequence (\ref{short1}) while the bottom one is induced by the short exact sequence 
\begin{equation}\label{eqlong}
0\rightarrow (\Omega^\bullet(M)[u], d+u\Delta) \overset{\times u}{\rightarrow} (\Omega^\bullet(M)[u], d+u\Delta) \overset{j}{\rightarrow} \Omega^\bullet(M)  \rightarrow 0,
\end{equation}
where $j(\sum a_i u^i)=a_0$, {\it cf.} \cite{GS,Ma}. In this picture, $m^\bullet$ corresponds to the connecting map of the long exact sequence \eqref{eqlong}. By a diagram chasing argument one finds that $m^\bullet=(i  \circ \Delta) ^\bullet$ where $i:\Omega^\bullet(M)\hookrightarrow \Omega^\bullet(M)[u]$ corresponds to $\mathcal B ^\bullet= (inc \circ B)^\bullet$  using Chen iterated integrals as corollary of Thoerem A in \cite{GJP}. Note that $i$  is not a chain map, whereas $i  \circ \Delta$ is a chain map, since $\Delta d=d\Delta$ and $\Delta^2=0$, ({\it cf.} \cite{GJP}).
\end{proof}

\section{$A_\infty$ generalization}\label{A-infty-section}
The previous sections, given for the case of dgas $(A,d,\cdot)$ with invariant inner product $\omega:A\otimes A\to k$, generalize in a straightforward way to the setting of cyclic $A_\infty$ algebras. In this section, we recall the relevant definitions  (\textit{cf.} \cite{T}), and adopt the above to this situation. 

\begin{defn}
An $A_\infty$ algebra on $A$ consists of a sequence of maps $\{\mu_n\}_{n\geq 1}$, where $\mu_n:A^{\otimes n}\to A$ is of degree $(2-n)$, such that
\begin{equation*}
\forall n\geq 1: \sum_{ \tiny\begin{matrix} {k+l=n+1}\\ r=0,\cdots,n-l \end{matrix} } (-1)^{\epsilon_{l}^r} \cdot \mu_k(a_1\otimes\cdots\otimes\mu_{l} (a_{r+1}\otimes \cdots\otimes a_{r+l})\otimes \cdots\otimes a_n )=0,
\end{equation*}
where $\epsilon_{l}^r=(l-1)\cdot (|a_1|+\cdots+|a_r|-r)$. A unit is an element $1\in k\subset A^0$ such that $\mu_2(a,1)=\mu_2(1,a)=a$, and $\mu_n(\cdots \otimes 1 \otimes \cdots )=0$ for $n \neq 2$. Again, we write $\bar A=A/k$. We define the Hochschild chain complex of $A$ with values in $A$ or $A^*$ to be the vector spaces $\bar C_\bullet(A,A)$ and $\bar C_\bullet(A,A^*)$ from equation \eqref{C(A,M)} with the differentials modified as follows,
\begin{eqnarray*}
\delta:\bar C_\bullet(A,A)\to\bar C_\bullet(A,A),  \delta(a_0\otimes\cdots \otimes a_n)=\sum \pm a_0\otimes\cdots\otimes \mu_k(\cdots)\otimes \cdots \otimes a_n\\
 + \sum\pm \mu_k(a_{s}\otimes \cdots \otimes a_0\otimes \cdots\otimes a_{r})\otimes a_{r+1}\otimes \cdots \otimes a_{s-1},\\
\delta:\bar C_\bullet(A,A^*)\to\bar C_\bullet(A,A^*), \delta(a_0^*\otimes\cdots \otimes a_n)=\sum \pm a_0^*\otimes\cdots\otimes \mu_k(\cdots)\otimes \cdots \otimes a_n\\
 + \sum\pm \mu_k^*(a_{s}\otimes \cdots \otimes a_0^*\otimes \cdots\otimes a_{r})\otimes a_{r+1}\otimes \cdots \otimes a_{s-1},
\end{eqnarray*}
where $ \mu_k^*(a_{s}\otimes \cdots \otimes a_0^*\otimes \cdots\otimes a_{r})\in A^*$ is given by $$\mu_k^*(a_{s}\otimes \cdots \otimes a_n\otimes a_0^*\otimes a_1\otimes \cdots \otimes a_{r}) (a):= \pm a^*_0 (\mu_k(a_1\otimes \cdots\otimes a_r\otimes a\otimes a_s\otimes \cdots \otimes a_n )). $$
Here, the signs are given by the usual Koszul rule, where we a factor of $(-1)^{\epsilon \epsilon'}$ is introduced, whenever elements of degree $\epsilon$ and $\epsilon'$ are being commuted. For an explicit discussion of the signs, see {\it e.g.} \cite{T}. Similarly, $\bar C^\bullet(A,A)$ and $\bar C^\bullet(A,A^*)$ are defined by the spaces from \eqref{C*(A,M)} with the modified differentials
\begin{multline*}
\delta^*:\bar C^\bullet(A,A)\to\bar C^\bullet(A,A), \quad \delta^* f(a_1\otimes\cdots \otimes a_n)\\=\sum \pm f(a_1\otimes\cdots\otimes \mu_k(\cdots)\otimes \cdots \otimes a_n)
 + \sum\pm \mu_k(a_1 \otimes \cdots\otimes f(\cdots)\otimes\cdots \otimes a_n),
\end{multline*}
\begin{multline*}
\delta^*:\bar C^\bullet(A,A^*)\to\bar C^\bullet(A,A^*), \quad \delta^* f(a_1\otimes\cdots \otimes a_n)\\
=\sum \pm f(a_1\otimes \cdots\otimes \mu_k(\cdots)\otimes \cdots \otimes a_n)
 + \sum\pm \mu_k^*(a_1\otimes \cdots\otimes f(\cdots)\otimes\cdots \otimes a_n).
\end{multline*}
Since $\delta^2=0, (\delta^*)^2=0$ in all the above cases, we obtain the associated homologies and cohomologies $H_\bullet(A,A), H_\bullet(A,A^*), H^\bullet(A,A)$, and $H^\bullet(A,A^*)$.

There is a generalized cup product $\cup$ on $H^\bullet(A,A)$ induced by,
$$ (f\cup g)(a_1\otimes \cdots\otimes a_{n}):=\sum_{k\geq 2} \pm \mu_k(a_1\otimes \cdots\otimes f(\cdots)\otimes\cdots\otimes g(\cdots)\otimes\cdots\otimes a_{n}). $$ Furthermore, equation \eqref{B} defines an operator $B:\bar C_\bullet(A,A)\to \bar C_\bullet(A,A)$ with $B^2=\delta B+B\delta =0$. We define the negative cyclic chains $\overline{CC}_\bullet^-(A)$ of $A$ to be  the vector space $\bar C_\bullet(A,A)[[u]]$ with differential $\delta+uB$, and denote the negative cyclic homology by $HC_\bullet^-(A)$. Dualizing $\overline{CC}_\bullet^-(A)$, we obtain  $\overline{CC}^\bullet_-(A)$ with dual differential and denote the negative cyclic cohomology by $HC^\bullet_-(A)$. For the same reasons as in Section \ref{sectionHC-}, we obtain the long exact sequences \eqref{CC_HH} and \eqref{CC^HH}.

Finally, assume we have a trace $\tr:A\to k$, such that the associated map $\omega:A\otimes A\to k, \omega(a,b)=\tr(\mu_2(a\otimes b))$ is a quasi-isomorphism, which satisfies for $n\geq 1$,
\begin{equation}\label{invariance}
\omega(\mu_n(a_1\otimes\cdots\otimes a_n),a_{n+1})=\pm \omega(\mu_n(a_{n+1}\otimes a_1\otimes\cdots\otimes a_{n-1}),a_n), 
\end{equation}
In this case, $\omega:A\to A^*$ induces a map of the Hochschild cohomologies $H^\bullet(A,A)\to H^\bullet(A,A^*), \omega^\bullet_\sharp(f)=\omega\circ f$, which we assume to be an isomorphism. Thus, we may transfer the product $\cup$ on $H^\bullet(A,A)$ to a product $\sqcup$ on $H^\bullet(A,A^*)$. $(HH^\bullet(A,A^*), \sqcup,\Delta=B^*)$ is a BV-algebra, {\it cf.} \cite{T}, so that we obtain the Lie bracket $\{a,b\}:=I^\bullet(\mathcal B^\bullet(a)\sqcup \mathcal B^\bullet(b))$ on $HC^\bullet_-(A)$ just as in Proposition \ref{Lie}.
\end{defn}

Using this setup, we may now also generalize Section \ref{MC-section}.
\begin{defn}
Recall that there are maps from the the $n^{th}$ symmetric power of a vector space to the $n^{th}$ tensor power  $S^n:A^{\wedge n}\to A^{\otimes n}$, where $S^n (a_1\wedge\cdots\wedge a_n)=\sum_{\sigma\in \Sigma_n} (-1)^{\epsilon_\sigma} (a_{\sigma(1)}\otimes \cdots\otimes a_{\sigma(n)})$. Defining $\nu_n:A^{\wedge n}\to A$ as $\nu_n:=\mu_n\circ S^n$, we obtain an $L_\infty$ algebra on $A$, {\it cf.} \cite[Theorem 3.1]{LM}. Furthermore, from \eqref{invariance}, it is immediate to see that we have for $n\geq 1$,
\begin{equation*}
\omega(\nu_n(a_1\wedge\cdots\wedge a_n),a_{n+1})=\pm \cdot\omega(\nu_n(a_{n+1}\wedge a_1\wedge \cdots\wedge a_{n-1}),a_n).
\end{equation*}

For this $L_\infty$ algebra, recall from \cite[Section 2]{GG} that the Maurer-Cartan solutions  are defined by,
\begin{eqnarray*}
 MC&:=& \Big\{a\in \A^{1}~\Big|~ \nu_1 (a)+\frac{1}{2!}\nu_2(a\wedge a) +\frac{1}{3!} \nu_3 (a\wedge a\wedge a)+\cdots =0\Big\}, \text{ and} \\
 \MC&:=& MC\big/\sim, 
\end{eqnarray*}
where the equivalence is again generated by the infinitesimal action of $\A^{0}$ on $\A^{1}$, where for $a\in A^{0}$, the vector field $\xi_{x}$ on $\A^{1}$ is defined by,
$$
\xi_{x}(a)=\nu_1(x) +\nu_2(a\wedge x)+\frac{1}{2!}\nu_3(a\wedge a\wedge x)+ \cdots.
$$
Note, that under the above assumptions the tangent space to $\MC$ at $[a]$ is the self-dual 3-term complex,
\begin{equation}
T_{[a]}\MC: \quad  T_{[a]}^{-1} \MC:=\A^{0}\overset{\xi(a)}{\longrightarrow}T_{[a]}^{0} \MC:=T_{a} A^{1}=A^{1}\overset{\mu'_{a}}{\longrightarrow} T_{[a]}^{1} \MC:={\A^{0}}^*,
\end{equation}
where 
$$\mu'_a(b)=\nu_1(b)+ \nu_2(a\wedge b)+\frac{1}{2!}\nu_3(a\wedge a\wedge b)+\cdots.$$ 
The self-duality at the middle term is given by the symplectic form
$$\omega(X_a,Y_a)=\tr(\mu_2(X_a\otimes Y_a))\in k.$$ 
This can be used to define the Hamiltonian vector field $X^\psi$ associated to a function $\psi\in \mathcal O(\MC)$, and thus the Lie bracket on $\mathcal O(\MC)$ via the usual formula $\{\psi, \chi\}=\omega(X^\psi,X^\chi)$.
\end{defn}
We may now define the map $P:MC\to \bar C_\bullet(A,A)$ by 
$$ P(a):=\sum_{i\geq 0}1\otimes a^{\otimes i}=(1\otimes 1_{\bar A^{\otimes 0}})+(1\otimes a)+(1\otimes a\otimes a)+\cdots, $$
and $R=inc\circ P:MC\to \overline{CC}^-_\bullet(A)$. As in definition \ref{PRrho}, we may again see, that $\delta(P(a))=0$, and $(\delta+uB)(R(a))=0$, and we define 
\begin{eqnarray*}
\rho:HC_-^{2\bullet} (\A)&\to& \mathcal O(\MC),\\
\rho([\alpha])([a])&:= &\langle [\alpha], [R(a)] \rangle=\langle \alpha, R(a) \rangle, \quad \quad \text {for } [\alpha]\in HC^\bullet (\A), [a]\in\MC.
\end{eqnarray*}
With this, we have the same theorem as in the previous sections.
\begin{thm}
The map $\rho$ is a well-defined map of Lie algebras.
\end{thm}

\address{Laboratoire Jean Leray, Universit\'e de Nantes, Nantes 44300, France.}
\email{abbaspour@univ-nantes.fr}

\address{College of Technology of the City University of New York, 11201 New York. Max-Planck Institut f\"ur Mathematik, Bonn 53111, Germany.}
\email{ttradler@citytech.cuny.edu}

\address{Long Island University, C.W. Post College, Brookville, NY 11548, USA.}
\email{mzeinalian@liu.edu}

\end{document}